\documentclass[11pt,a4paper,leqno]{amsart}
\usepackage{amssymb, amscd}
\usepackage[all, cmtip]{xy}
\newtheorem{theorem}{Theorem}[section]

\newtheorem{corollary}[theorem]{Corollary}
\newtheorem{lemma}[theorem]{Lemma}
\newtheorem{definition}[theorem]{Definition}

{\sc}

{}
{}
{\sc}
\newtheorem{remark}{Remark}[section]

\newcommand{\acknowledge}{\subsection*{Acknowledgments}}
\newcommand{\thismonth}{\ifcase\month\or
  January\or February\or March\or April\or May\or June\or
  July\or August\or September\or October\or November\or December\fi
  \space\number\year}

\makeatletter
\newcommand{\rssymb}[2]{\newcommand{#1}{{\mathrmsl{#2}}}}
\newcommand{\calsymb}[2]{\newcommand{#1}{{\mathcal{#2}}}}
\newcommand{\bbsymb}[2]{\newcommand{#1}{{\mathbb{#2}}}}
\newcommand{\lieoper}[2]{\newcommand{#1}{\mathop
  {\mathfrak{#2}\null}\nolimits}}
\newcommand{\oper}[3][n]{\newcommand{#2}{\mathop
  {\mathrm{#3}\null}\ifx n#1\nolimits\else\limits\fi}}
\newcommand{\rsoper}[3][n]{\newcommand{#2}{\mathop
  {\mathrmsl{#3}\null}\ifx n#1\nolimits\else\limits\fi}}
\bbsymb\C{C} \bbsymb\F{F} \bbsymb\HQ{H}\bbsymb\N{N} \bbsymb\Q{Q}
\bbsymb\R{R} \bbsymb\U{U} \bbsymb\V{V} \bbsymb\W{W} \bbsymb\Z{Z}
\bbsymb\bbf{F} \bbsymb\bbk{K} \bbsymb\bbi{I} \bbsymb\bbl{L}
\bbsymb\bbo{O} \bbsymb\bbj{J} \bbsymb\bby{Y} \bbsymb\bbp{P}
\bbsymb\bba{A}
\calsymb\cA{A} \calsymb\cB{B} \calsymb\cC{C} 
\calsymb\cM{M} \calsymb\cN{N} \calsymb\cO{O} \calsymb\cP{P}
\calsymb\cU{U} \calsymb\cV{V} \calsymb\cW{W} \calsymb\cX{X}
\calsymb\cY{Y} \calsymb\cZ{Z}

\newcommand{\lam}{\lambda}
 \renewcommand{\leq}{\leqslant}
\oper\End{End} \oper\Hom{Hom}                    
\oper\Sym{Sym} \oper\Skew{Skew}
\oper\Aut{Aut}                                   
\oper\GL{GL} \oper\SL{SL}\oper\Symp{Sp} \oper\CO{CO} \oper\On{O}
\oper\SO{SO} \oper\Pin{Pin} \oper\Spin{Spin} \oper\CU{CU}
\oper\Un{U} \oper\SU{SU} \oper\PSU{PSU} \rsoper\Diff{Diff}
\rsoper\SDiff{SDiff}
\lieoper\der{der}                                
\lieoper\gl{gl} \lieoper\sgl{sl}\lieoper\symp{sp} \lieoper\co{co}
\lieoper\so{so} \lieoper\spin{spin} \lieoper\cu{cu} \lieoper\un{u}
\lieoper\su{su} \rsoper\Vect{Vect} \rsoper\Ham{Ham}
\def\la#1{\hbox to #1pc{\leftarrowfill}}
\def\ra#1{\hbox to #1pc{\rightarrowfill}}

\newcommand{\norm}[2][]{|\mkern-2mu|#2|\mkern-2mu|
  _{\lower1pt\hbox{${}_{#1}$}}}
\newcommand{\Norm}[2][]{\bigl|\mkern-3mu\bigr|#2\bigr|\mkern-3mu\bigr|
  _{\lower1pt\hbox{${}_{#1}$}}}

\rsoper\dimn{dim}                           
\rsoper\grad{grad}                          
\rsoper\kernel{ker}\rsoper\image{im}        
\rsoper\alt{alt}   \rsoper\sym{sym}         
\rsoper\Ad{Ad}     \rsoper\ad{ad}           
\rsoper\CoAd{CoAd} \rsoper\coad{coad}       
\rsoper\trace{tr}  \rsoper\trfree{tf}       
\rsoper\detm{det}                           
\rsoper\Vol{Vol}                            
\rsoper\divg{div}                           
\rsoper\sign{sign}                          
\rssymb\iden{id}                            
\rssymb\vol{vol}                            
\oper\Imag{Im}\oper\Real{Re}                
\newcommand{\sd}{{\raise1pt\hbox{$\scriptscriptstyle +$}}}
\newcommand{\asd}{{\raise1pt\hbox{$\scriptscriptstyle -$}}}
\newcommand{\sdasd}{{\raise1pt\hbox{$\scriptscriptstyle\pm$}}}
\newcommand{\asdsd}{{\raise1pt\hbox{$\scriptscriptstyle\mp$}}}

\rsoper\scal{scal}
\def\kahl/{k\"ahler}
\def\Kahl/{K{\"a}hler}

\begin{document}
\title[Null Sasaki $\eta$-Einstein Structures in Five Manifolds]
{Null Sasaki $\eta$-Einstein Structures in Five Manifolds}
\author[J. Cuadros]{Jaime Cuadros Valle}
\address{Department of Mathematics and Statistics,
McMaster University,
Hamilton, On L8S 4K1, Canada.}
\email{jcuadros@math.mcmaster.ca}

\date{\thismonth}
%
%

\begin{abstract} We study null Sasakian structures in dimension five. First, based on a result due to Koll\'ar \cite{Ko}, we improve a result in \cite{BGM} and prove that simply connected manifolds diffeomorphic to $\# k(S^2\times S^3)$ admit null Sasaki $\eta$-Einstein structures if and only if $k\in \{3,\ldots , 21\}$. 
After this, we determine the moduli space of simply connected null Sasaki $\eta$-Einstein structures. This is accomplished using information on the moduli of lattice polarized K3 surfaces.
\end{abstract} 
\maketitle
\vspace{-2mm} 




\section{Introduction}

The purpose of this paper is to classify {\it null Sasaki $\eta$-Einstein} five manifolds. Recall that to every Sasakian manifold corresponds 
a basic first Chern class
$c_1({\mathcal F}_{\xi})$ which leads to define null, positive or negative Sasakian structures 
in analogy with K\"ahler geometry. 
As a consequence of a transverse version of Yau's celebrated theorem \cite{EKA} every null Sasakian structure can be deformed to a Sasaki 
$\eta$-Einstein structure which is transverse Calabi-Yau. One refers to these structures as null Sasaki $\eta$-Einstein.   

In \cite{BGM}, the idea of links $L_f$ associated to Reid's list of all families of codimension one K$3$ 
surfaces in weighted projective spaces $\mathbb{P}(w_0,w_1,w_2,w_3)$ is used to exhibit several examples of null Sasaki $\eta$-Einstein structures. However, we need to extend this idea to links associated to complete intersections to give a complete classification of these structures. For this purpose one needs to determine the second Betti number of these links  $L_{\bf f}$. One way of doing this for links coming from weighted hypersurfaces is using the Milnor-Orlik formula  \cite{BGK}. But this method is not applicable to links coming from weighted complete intersections. 
We find a way to calculate these Betti numbers: one considers a  Seifert bundle as a 
``resolution'' of an orbifold. How far is this odd dimensional resolution from the usual minimal resolution. 
When the orbifold fundamental group is trivial and under the absence of branch divisors these two notions are, in certain sense, very close. Under this setting, we exploit some well-known  properties of Du Val singularities (the ones that appear into play) 
to determine these Betti numbers.  
We appeal to a  Smale's classification theorem on simply connected 
spin five dimensional manifolds \cite{Sm} and a result due to  Koll\'ar \cite{Ko} on the absence of brach divisors to  improve previous results on the  classification of null Sasaki $\eta$-Einstein five manifolds given in \cite{BGM}. In  Section 3 we prove 

\noindent{\bf Theorem 1. }{\it 
 Let $M$ be a simply connected five dimensional manifold that admits a null Sasaki $\eta$-Einstein structure.
 Then $M$ is diffeomorphic to a connected sum of $k$ copies of $S^2\times S^3$, where $k$ is the second Betti 
 number of $M$ and $k$ takes every integer between 3 and  21.} 

In Section 4 we describe the space of deformations of null Sasaki $\eta$-Einstein five manifolds.  
For this we deform the transverse holomorphic Sasakian structures. These deformations arise 
from deforming the complex structures on the base space keeping a fixed 
(locally free) action of the circle on $M$ and therefore can be seen as deformations of the orbifold keeping the orbifold singularities, 
produced by this action,  fixed. This has an interpretation in terms of deformations of lattice polarized K3 surface: 
we appeal to the deformations of the minimal resolutions of each of the 95 members of the family on Reid's list and 84 members of the 
family of Fletcher's list of weighted complete intersections keeping the corresponding exceptional loci fixed (or polarized). We also need to consider an analog to the quadratic form in K3 surfaces in terms of basic 2-forms. 
Thanks to a result of Tondeur \cite{Ton}, this quadratic form is non-degenerate. It is possible to split the representative in cohomology in terms of a class in $H^2(M, \mathbb C)$ and the fundamental class $[d\eta]_B\in H^2({\mathcal F}_\xi)$ and in this way we  find some information on the {\it inherited} 
quadratic form in $H^2(M, \mathbb C)$. We obtain

\noindent{\bf Theorem 2.} {\it  
Simply connected five dimensional manifolds $M$ that admit null Sasaki $\eta$-Einstein metrics have a space of deformation of dimension equal to $2(b_2(M)-2)$. Moreover, the moduli space of null Sasaki $\eta$-Einstein structures for $M$ is determined by the following quadric in $\mathbb C\mathbb P^{b_2(M)-1}$: 
$$\{ [\alpha]\in H^2(M, \mathbb C)\, |  \, \rVert([\alpha], [\alpha])\rVert-\rVert([\alpha], [ \overline{\alpha}])\rVert >0\}/\mathbb C^*.$$} 


In Section 2 we go over some preliminaries on complex orbifolds. 
We also recall the definition of a Sasakian manifold, invariants of 
its characteristic foliation and some properties on links of complete intersections. 

\acknowledge Some of these results were part of the author's dissertation 
under the supervision of Prof. Charles P. Boyer. I would like to thank him and Krzysztof Galicki for their invaluable help. 





\section{Preliminaries}

In this section we review some basic definitions on orbifolds, Sasakian structures and their relations, see details in \cite{BBG}.

\subsection{Basics on Orbifolds.} Recall that on an orbifold $(X, \Delta)$, the branch divisor 
$\Delta=\sum_i(1-\frac{1}{m_i})$ is a $\mathbb Q$ divisor on the underlying normal compact complex space $X$. The $D_i$'s  are Weil divisors and the ramification index $m_i$, is the greatest common divisor of the the orders of all the uniformizing groups of $X$ at every point along $D_i$. 
The very important orbifold canonical divisor $K_X^{orb}$  is related  
to the canonical divisor of the underlying complex structure by the formula 
\begin{equation}
K^{orb}_X=K_X+\Delta, 
\end{equation}
hence the orbifold Chern class is given by $c_1^{orb}(X, \Delta):=c_1(X)-\Delta.$ When the fixed point set has codimension at least 2, then $\Delta\equiv 0$ and we have  $K^{orb}_X=K_X.$

Another important topological invariant of an orbifold is given by the {\it orbifold fundamental group} 
$\displaystyle{\pi_1^{orb}(X, \Delta)}$, is the fundamental group of $X\backslash ({\mathrm {Sing}} X\cup {\mathrm {Supp}} \Delta)$ 
modulo the relations: if $\gamma$ is any small loop around $D_i$ at a smooth point then $\gamma^{m_i}=1$. 
The abelianization of $\pi_1^{orb}(X, \Delta)$ is denoted by $H_1^{orb}(X, \Delta)$. 
\medskip


\noindent{\it Complete Intersections in Weighted Projective Spaces ${\mathbb P}$(\bf w).}
Interesting examples of K\"ahler orbifolds --and ultimately of Sasakian structures in odd dimensional manifolds-- 
are given by complete intersections in $\mathbb P({\bf w})$. 
A {\it weighted variety} X is the common zero locus in $\mathbb P({\bf w})$ of a collection of weighted homogeneous polynomial in 
$\C[z_0, \ldots z_n]$. A {\it weighted hypersurface} $X_f$ is the zero locus in  $\mathbb P({\bf w})$ of a single weighted homogeneous 
polynomial $f$. A weighted variety is called {\it weighted complete intersection}  if the number of polynomials in the collection equals the 
codimension of $X$.  
$X_d\subset {\mathbb P}({\bf w})$ will denote the weighted hypersurface of degree $d$ and with 
$X_{d_1, \ldots , d_c}\subset{\mathbb P}({\bf w})$ the weighted complete intersection of multidegree ${d_i}$ 
(here $c$ denotes the codimension of $X$).  
{\it Quasismooth} weighted complete intersections are those with at worst cyclic singularities and, thus, can be seen as orbifolds 
(see \cite{IF} for details). In this situation,  the well-formedness of $X$ implies that 
the orbifold singular locus and the algebro-geometric singular locus coincide. Later, in proving Theorem 1 and 2 we will use weighted varieties 
with these properties. 


\subsection{Sasaki $\eta$-Einstein Structures.} 
A $2m+1$ dimensional  manifold $(M,g)$ is said to be {\it Sasakian} if its corresponding cone $C(M)={\mathbb R}_+\times M$ with the warped 
metric $\overline{g}=dr^2+r^2g$ is K\"ahler. From this definition 
we obtain a vector field $\xi$ (the Reeb vector field) and a one form 
$\eta$ (the contact form)  both on $M$. The one dimensional foliation generated by $\xi$ is called the {\it Reeb foliation} and denoted by 
$\mathcal{F}_\xi$. It is possible to extend this vector field to  a holomorphic vector field on $C(M)$, thus there is a holomorphic action of 
${\mathbb C}^*$ on $C(M)$. The local orbits of this action define a transversely holomorphic structure on the Reeb foliation. 
Furthermore, the Reeb foliation is transversely K\"ahler. 
There is also a $(1,1)$ tensor $\Phi$ that can be defined in terms of the complex structure in the cone. 
Indeed, the complex structure of the transversely K\"ahler structure is given by $\Phi$ restricted to the kernel of $\eta$. 
We denote by $(M,\xi, \eta, \Phi, g)$ the Sasakian manifold.  


An important transverse invariant is the {\it basic first Chern class} $c_1({\mathcal F}_\xi)\in H^2_B({\mathcal F}_\xi)$ of the foliation ${\mathcal F}_\xi$.  Here  $H^*_B({\mathcal F}_\xi)$ denotes the {\it basic cohomology ring} of smooth basic forms of a foliation  $\mathcal{F}_\xi$ (and since the transverse geometry is K\"ahler, we also have transverse Hodge decompositon at our disposal). 

\begin{definition}
A Sasakian structure $(\xi,\eta, \Phi, g)$ is said to be null if $c_1({\mathcal F}_\xi)=0,$  that is the transverse structure is Calabi-Yau. 
It is {\it positive (negative)} if $c_1({\mathcal F}_\xi)$ is represented by a positive (negative) definite $(1,1)$  form. 
If either of these two conditions is satisfied $(\xi, \eta, \Phi, g)$ is said to be {\it definite}. 
\end{definition}

The next definition, introduced in \cite{Ok}, gives us an anolog to K\"ahler-Einstein structures for odd dimensions. 
\begin{definition}
A Sasakian structure ${\mathcal S}=(\xi, \eta, \Phi, g)$ on $M$ is said to be {\it Sasaki $\eta$-Einstein} or just 
$\eta$-Einstein if there are constants $\lambda, \mu$ such that $\mathrm{Ric}_g=\lambda g+\mu\eta\otimes \eta.$ 
It follows that every Sasaki 
$\eta$-Einstein manifold is of constant scalar curvature $s=2n(\lam+1).$
\end{definition}

As in the K\"ahlerian case,  we have  that the basic class $2\pi c_1({\mathcal F}_\xi)\in H_B^{1,1}({\mathcal F}_\xi)$ is represented 
by the transverse Ricci form $\rho_g^T.$ The converse to this last statement is also true. 
In complete analogy with the K\"ahlerian case, we see that there are no obstructions to solving the Sasakian 
or the transverse Monge-Amp\`ere problem in the negative or null case, see \cite{EKA} for details. 
\medskip

\noindent{\it Null $\eta$-Einstein Structures on Links.} 
Now we collect some information about links of complete intersections defined by $p$ linearly indepedent weighted homogeneous polynomials 
$f_1, \ldots , f_p \in \C^{n+1}$ of degrees $d_1, \ldots , d_p$, respectively with weight vector ${\mathbf w}$. 
Consider the {\it weighted affine cone} $C_{\mathbf f}=\{ (z_0, \ldots , z_n)\in \C^{n+1} |\,  {\mathbf f}(z_0, \ldots , z_n)=0\},$ 
which has dimension $n+1-p$. Let us assume  that the origin in $\C^{n+1}$ is the only singularity, and that this singularity is isolated. 
Then we define  the link $ L_{\mathbf f}=C_{\mathbf f}\cap S^{2n+1},$ which is smooth, of dimension $2(n-p)+1$ and it was proven in \cite{Lo} that it is $(n-p-1)$-connected. 
If ${\mathcal S}_{\bf w}=(S_{\bf w}^{2n+1},\xi_{\bf w}, \eta_{\bf w}, \Phi_{\bf w}, g_{\bf w})$ denotes 
the natural weighted Sasakian structure on the sphere corresponding to  ${\mathbb P}({\bf w})$,  
then the quadruple $(\xi_{\bf w}, \eta_{\bf w}, \Phi_{\bf w}, g_{\bf w})$  endows  $L_{\bf f}$ with a quasi-regular  
Sasakian structure. We have the following commutative diagram
\begin{equation}
\begin{CD}
L_{\mathbf f}  @> {\qquad\qquad}>>  S^{2n+1}_{\bf w}\\
@VV{\pi}V  @VVV\\
{\mathcal Z}_{\mathbf f}  @> {\qquad\qquad}>>  {\mathbb P}({\bf w})
\end{CD}
\end{equation}
Here the horizontal arrows are Sasakian and K\"ahlerian embeddings, respectively. 
Moreover $\pi:L_{\mathbf f}\rightarrow {\mathcal Z}_{\mathbf f},$ is an orbifold riemannian 
submersion and $\pi^*c_1^{orb}(\mathcal{Z}_{\bf f})=c_1({\mathcal F}\xi).$
The algebraic variety ${\mathcal Z}_{\mathbf f}$ is a weighted complete intersection in the weighted projective space 
${\mathbb P}({\bf w})$, and the condition that the cone $C_{\mathbf f}$ be smooth outside the origin is translated into the 
quasi-smoothness of ${\mathcal Z}_{\mathbf f}$.
Let us denote by  ${\mathcal S}_{{\bf w}, {\mathbf f}}=S_{\bf w}|_{L_{\mathbf f}}$ 
the Sasakian structure endowed from ${\mathcal S}_{{\bf w}}$ on the link $L_{\mathbf f}$. From the orbifold adjunction formula 
$K^{orb}_{\mathcal{Z}_{\bf f}}=\mathcal{O}_{\mathcal{Z}_{\bf f}}(\sum d_i -\sum w_i)$, it follows that ${\mathcal S}_{{\bf w}, {\mathbf f}}$ 
is null when $\mathcal{Z}_{\bf f}$ is null,  that is, when $\sum d_i -\sum w_i=0$. Moreover, the lack of obstructions in the non-positive case 
implies that $L_{\mathbf f}$ admits a Sasakian structure 
$\mathcal{S}'$ with $\eta$-Einstein metric $g'$ in the same deformation class as ${\mathcal S}_{{\bf w}, {\mathbf f}}.$ 

\section{Seifert Bundles and K3 Surfaces}

A compact K\"ahler orbifold $(\cZ, \omega)$ is called 
{\it Hodge orbifold} if $[\omega]$ lies in $H_{orb}^2(\cZ, \Z).$ By a theorem of Baily \cite{Ba}  a Hodge orbifold is a polarized 
projective algebraic variety. The {\it K\"ahler lattice} $\mathcal{K}_L(\mathcal{Z})$ is defined as the intersection of the K\"ahler cone with  
$H_{orb}^2(\mathcal{Z}, \mathbb Z)$. 
Boyer and Galicki proved
\begin{theorem}[\cite{BBG}]
Let $(\cZ, \omega)$ be a Hodge orbifold. Every $[\omega]\in\mathcal{K}_L(\mathcal{Z})$ determines a Seifert $S^1$-bundle or equivalently a 
principal $S^1$ V-bundle $\pi: M_{[\omega]}\longrightarrow \cZ$ and choosing a connection $\eta$ in $M_{[\omega]}$ whose curvature is 
$\pi^*\omega$ determines a Sasakian structure on $M_{[\omega]}.$ Moreover, if the orbifold is locally cyclic 
(that is,  if it has an orbifold atlas all whose local uniformizing groups are cyclic groups) then $M$ is a manifold.
 \hfill$\blacksquare$
\end{theorem}

The following theorem will be used to prove Theorem 1. Its proof is a direct consequence of a classification theorem on  simply connected spin five manifols due to  Smale in \cite{Sm} and  Theorem  5.7 in \cite{Ko}. 

\begin{theorem}
Let $f:M \longrightarrow (X, \Delta)$ be a simply connected five dimensional smooth Seifert bundle over a projective surface with rational singularities. Then the cohomology groups $H^i(M, \Z)$ 
are explicitly described. In particular, under the absence of branch divisors we have 
$$
H^2(M,\Z)=H^3(M,\Z)=\Z^{b_2(X)-1}.   
$$
Moreover, if $M$ is spin, it is diffeomorphic to the connected sum of $b_2(X)-1$  copies of $S^2\times S^3.$
\hfill$\blacksquare$
\end{theorem}

\subsection{K3 Surfaces}
We have to display some facts about K3 surfaces (excellent reference is  \cite{BPV}). 
Recall 
that a  ${\mathrm K3}$ surface $X$ is a compact K\"ahler surface with only Du Val singularties or rational double points  such that $h^1({\mathcal O}_X)=0$  and the dualizing sheaf $\omega_X$ is trivial, that is $\omega_X={\mathcal O_X}$. 
If $X$ is non-singular, the dualizing sheaf becomes the line bundle associated to the canonical line bundle $K_X$. In this case, $H^2(X, \Z)$ is torsion free of rank 22.  
By means of intersection form, $H_2(X, \Z)$ is endowed with the structure of a lattice. 
It is isomorphic to the {\it K3 lattice} $L_\Lambda=-E_8^2\oplus H(1)^3,$ 
where $H(m)$ is the indefinite rank 2 lattice with intersection {\tiny{$\left ( \begin{array}{cc}
{0} & {m} \\
{m} & {0}
\end{array} \right )$} } and $E_8$ is the root-lattice associated to the Dynkin diagram $E_8$ 
which is positive definite. By Poincar\'e duality, $H^2(X,\Z)$ is also equipped 
with the lattice structure via the cup-product. This product is even, unimodular and  indefinite with signature (3,19). 
From the Hodge index theorem, the signature of the NS$(X)=H^{1,1}(X, \R)\cap H^2(X, \Z)$ is $(1, \rho-1),$  
where $\rho={\mathrm{rank}}NS(X)$. 
Non-singular rational curves can be blown down to yield rational double points: they are one of the following type $A_n, D_n, E_6, E_7, E_8$.  Resolution of these singularities have exceptional locus with 
transversally intersected smooth rational curves, forming a basis for a subspace included in $H^2(X, \mathbb{Z})$. 
From the adjunction formula we know that a smooth curve $C$ on a non-singular K3 surface is rational if and only if $C.C=-2$. 
Moreover, any irreducible curve has self-intersection 0 (mod 2). 

A classical example of a smooth  K3 surface is a {\it Kummer Surface} defined as the minimal resolution of the 16  singularities of 
type $A_1$ of $Z/\iota$, the quotient of a  complex torus $Z$ of dimension 2 by an involution $\iota$ on $Z$ 
(induced by multiplication by $-1$ on ${\mathbb C}^2$). 
\medskip

\noindent{\it Weighted K3 Surfaces of Codimension One and Two.} 
As explained before, links of complete intersections in dimension five admitting Sasakian structures will have as base space a 
weighted surface complete intersection which naturally has the structure of a K\"ahler orbifold. 
In the last pages of this work we give the lists of well-formed quasi-smooth weighted surface complete intersections 
of codimension 1 and 2 found by Reid and Iano-Fletcher respectively. All the members of these families satisfy 
$\sum_i d_i-\sum_j a_j=0$ and hence are K3 surfaces with at worst Du Val singularities. 
We also give the second Betti number for each member of these families. As we will see, 
these two lists ensure the existence of $\#k(S^2\times S^3)$  
admitting null Sasaki $\eta$-Einstein for every $3\leq k\leq 21$. 
The following lemma will be used in the proof of Theorem 1. 
\begin{lemma}
There are no projective K3 surfaces with second Betti number 3 and only rational double points. 
\end{lemma}
\noindent{\bf Proof.} 
Let us suppose there is a K3 surface with the features given above. Then its minimal resolution $X$ must have a  $NS(X)$ with configuration 
consisting of nineteen  $-2$ curves coming from the singularities, and of the hyperplane bundle. Hence $NS(X)$ is maximal, that is, its rank 
is 20, with signature (1,19), and it is even. 
Considering now its transcendental lattice $T_X$, with signature $(2,0)$ which in this case is also even. Denote by  $v_1, w_1, v_2, w_2$ 
with $(v_1, v_2)=(w_1,w_2)=0$ and $(v_i,w_j)=\delta_{ij}$ a basis of $H\oplus H$. If $\{ l_1,l_2 \}$ is a basis of $T_X$ then we define $f$ 
via the rule
$$f(l_1)=v_1+\frac{1}{2}(l_1, l_1)w_1 \quad \hbox{and }\, f(l_2)=v_2+\frac{1}{2}(l_2, l_2)w_2+ 2(l_1, l_2)w_1$$ 
which satisfies $(f(l_1), f(l_2))=2(l_1,l_2)$. From the fact that $(f(l_i), w_j)=\delta_{ij}$ we conclude that $f$ is a primitive embedding. 
It follows that $T_X$ has a primitive embedding in $H(2)^3$. 
This implies that  its complement $(H(2)^3)^{\perp}$ in $L_{\Lambda}$ is a primitive embedding in $NS(X)$. Notice that $(H(2)^3)^{\perp}$ is even, 
negative-definite and of rank 16. In \cite{Ni} Nikulin proved that this is the case if and only if $X$ is a Kummer surface. 
This is not the situation, otherwise the projective torus that gives rise to the Kummer surface has 3 more rational 
curves (we need 19 of these curves). However, it is known that abelian varieties admit no  rational curves (see \cite{BL}, Proposition 4.9.5). 
\hfill$\blacksquare$

\medskip


\subsection{Null Sasaki $\eta$-Einstein Structures in Dimension Five} Here we specialize to the case of interest in this article, 
that is Seifert bundles with base space a Calabi-Yau orbifold with at worst cyclic quotient singularities. 

\begin{definition}
 A compact orbifold $(X,\Delta)$ that admits a K\"ahler metric is called a Calabi-Yau orbifold if $\pi_1^{orb}(X)=1$ and its orbifold canonical class $K_X^{orb}$ is (numerically) trivial. 
\end{definition}

The next theorem ensures the absence of branch divisors under mild conditions.  
\begin{theorem}[\cite{Ko}]
 Let $f:M\longrightarrow (X,\Delta)$ be a Seifert bundle with $M$ a smooth manifold and $(X, \Delta)$ a Calabi-Yau orbifold. If $H_1(M, \mathbb Z)=0$ then $\Delta=\sum_{i=1}(1-\frac{1}{m_i})=0$. Thus, the canonical class $K_X^{orb}=K_X$ is trivial. 
\hfill$\blacksquare$
\end{theorem}

Next, we present the proof of Theorem 1 which is an improvement to the original idea given in \cite{Ko} 
and \cite{BGM}: this new approach, ensures the existence of null Sasaki $\eta$-Einstein structures for $k=17$ 
and discards the possibility for $k=2$.  The key here (and also when proving Theorem 2) is to consider a  Seifert bundle as an odd dimensional  
resolution of the orbifold. How far is this resolution from the usual resolution (in the sense of algebraic geometry)?  
When the orbifold is well-formed and simply connected these two notions are, in certain sense, 
very close one to the other: one can ``identify'' connected sums of $S^2\times S^3$ with exceptional curves of the minimal resolution. 
 
\medskip

\noindent{\bf Theorem 1.} 
{\it Let $f:L\longrightarrow (X,\Delta)$ be a Seifert bundle with $L$ a smooth simply connected five dimensional manifold and 
$(X,\Delta)$ a Calabi-Yau orbifold. Then $L$ is diffeomorphic to $\#k(S^2\times S^3)$ for every integer $k$ between 3 and 21.} 
\medskip

\noindent{\bf Proof.}  
In \cite{HS} it is shown that a  Seifert bundle with $L$ smooth gives the exact sequence 
\begin{equation}
\pi_2^{orb}(X,\Delta)\longrightarrow \Z \longrightarrow \pi_1(L)\longrightarrow \pi^{orb}_1(X, \Delta) \longrightarrow 1, 
\end{equation}
since $L$ is simply connected, it follows that $X$ is necessarily simply connected, in the orbifold sense. 
From Theorem 3.7 we have the triviality  of $K_X$. 
In particular, the (cyclic quotient) singularities also have trivial canonical class and thus are Du Val. 
Let us consider $\rho:\widetilde{X}\longrightarrow X$, the minimal resolution of $X$; this resolution is crepant, 
{\it i.e.,} $K_{\widetilde X}=\rho^*K_X$. Also the  fundamental group remains unchanged, hence 
$\pi_1(\widetilde{X})=\pi_1^{orb}(X)=\pi_1(X\backslash Sing(X))=1$. From Kodaira's classification of compact complex surfaces, 
$\widetilde{X}$ is a K3 surface.

From Theorem 3.6 it follows the absence of branch divisors. Moreover, in \cite{Mo} it is shown that the existence of Sasaki $\eta$-Einstein structures  forces the second Stiefel-Whitney class 
to vanish. Applying Theorem 3.2 we conclude that $H_2(L, \Z)$ is torsion free and $b_2(X)=b_2(L)+1.$
The resolved singularities create smooth rational curves of self-intersection $-2$, 
those only  decrease the second Betti number on $\widetilde{X}$ and thus $\mathrm{rank} H^2(X, \Q)\leq \mathrm{rank} H^2(\widetilde{X}, \Q)$.  Again by Theorem 3.2, $L$ is diffeomorphic to a connected sum of at most $21$ copies 
of $S^2\times S^3$. 
To prove that $k$ takes every integer between 3 and 21, it is enough to show that there are four dimensional Calabi-Yau orbifolds
$X_{{\bf w}_i}$ 
(with no branch divisor) such that the corresponding total spaces admit  second Betti numbers ranging from 3 to 21 inclusive. 
This is equivalent to finding four dimensional Calabi-Yau orbifolds with second Betti numbers between 4 to 22.  
Our source of examples will be given by Reid's list of codimension 1 weighted K3 surfaces, 95 in all, and Fletcher's list of 84 codimension 2 
weighted K3 surfaces \cite{IF} (see tables 1 and 2). These K3 surfaces have no branch divisors, moreover their singularities are at worst  cyclic quotient canonical surface 
singularities, that is, they are of type $A_{n_i}$. Let $f:X\rightarrow X_{\bf w}$ be a minimal resolution. Then $X$ is a non-singular K3 
surface, hence $b_2(X)=22$. Since Du Val singularities contribute exactly with $\sum_{i}n_i$ lineraly independent curves on $X$ with self-intersection $-2$, 
it follows that $b_2(X_{\bf w})=22-\sum_{i}n_i$. As shown on the Table 1,  Reid's list provides singularities that give us  $b_2(X_{\bf w})$ 
ranging from 4 to 22 inclusive except 18. We appeal to Fletcher's list to ensure the existence of a weighted K3 surface --this time of 
codimension 2-- with second Betti number 18 (numbers 5, 6, 7, 9, 11 and 17 on Table 2 are the  orbifolds having  that property). 
As a consequence of Lemma 3.5, the possibility $b_2(L)=2$ is discarded.  
\hfill$\blacksquare$

\begin{remark}
{\rm It is important to notice that the Milnor-Orlik formula is not useful to calculate second Betti numbers of links associated, 
via Seifert bundle, to members of Fletcher's list.} 
\end{remark}

From \cite{BGK}, \cite{Go} and the previous theorem we have the following corollary. 

\begin{corollary}
Suppose a simply connected five-manifold $M^5$  admits  $\eta$-Einstein metrics of positive, negative, and null type. Then $M^5=\#k(S^2\times S^3)$ if and only if  $k\in \{3,\cdots,21\}.$
\hfill$\blacksquare$
\end{corollary}

But, how many inequivalent null Sasaki-Seifert structures  can one find in 
$\#k(S^2\times S^3)$ in each $k$?  For example, for $k=21$, where the action is regular, $\#21(S^2\times S^3)$ is the circle bundle 
over a non-singular projective K3 surface. For instance, numbers 1 and 3 in Table 1, number 1 in Table 2, and 
the intersection of three quadrics in $\mathbb P^5$ have these properties. But of course these are not the only ones. A complete answer will require a complete 
classification of simply connected projective K3 orbifolds (here the term K3 orbifold refers to those whose minimal resolution is a smooth 
K3 surface) with cyclic singularities, which at this moment is unknown to the author. (However, see \cite{KeZh} for a partial classification for 
$K3$ surfaces with finite fundamental group.)  
\medskip

\noindent{\it The Non-Simply Connected Case.} 
Theorem 3.5 establishes the rigidity of simply connected smooth null Sasakian manifolds. Its corresponding base space, in general an orbifold, 
has no branch divisor. However, a (projective) Calabi-Yau orbifold with singularities of codimension 1  will produce a Sasaki-Seifert manifold $L$ 
with null $\eta$-Einstein metric, but this time our manifold is not simply connected. Using the exact sequence (3) 
--and from the fact that $H_1(L,\mathbb Z)$ cannot be zero--, we have that $\pi_1(L)$  is either ${\mathbb Z}_k$ or $\mathbb Z$. 

In general, for a normal surface with no branch divisor and whose minimal resolution is a K3 surface is known (see \cite{KeZh} or \cite{Cam}) 
that its orbifold fundamental group (that is, the fundamental group of its smooth locus) is either finite or an extension of 
${\mathbb Z}^4$. 
An  example of the latter is given by the Kummer surface $X$ which has orbifold fundamental group spanned by ${\mathbb Z}^4$ 
and ${\mathbb Z}_{2}$ and null (orbifold) first Chern class. Clearly $X$ has 
no codimension 1 singularities. Notice that in this case the usual fundamental group is trivial, so the Seifert bundle viewed as resolution of 
the orbifold is not as close as desired to the the minimal resolution (a smooth K3 surface). If we apply sequence (3) 
we see that the corresponding 
null Sasaki $\eta$-Einstein manifold for $X$ (the uniformizing groups inject into $S^1$) cannot be a connected sum of $S^2\times S^3$. 
An explicit way to construct these structures in the smooth case arises from  compact Heisenberg manifolds 
(which are the quotient of the Heisenberg group by a lattice subgroup of it). It is not difficult to see that these manifolds are 
circle bundles  over an abelian variety. Folland proved that these Heisenberg manifolds are in one to one correspondence, 
up to holomorphic equivalence, with polarized abelian varieties \cite{Fol}. 
\section{Deformations of Transverse Holomorphic Sasakian Structures}

The cohomological groups with respect to the sheaf of holomorphic vector fields $T_X$ (on a complex manifold $X$) have precise interpretation 
in terms of deformation theory \cite{Kod}: 
$H^1(X, T_X)$ can be seen as  the space of infinitesimal deformations of the complex structure of $X$, $H^2(X, T_X)$ 
is the space of obstructions to lifting an infinitesimal deformation to an actual deformation of the complex structure of $X$.  
The group $H^0(X, T_X)$ can be considered as the Lie algebra of the group of holomorphic automorphisms of $X$, and if $H^0(X, T_X)=0$ 
then the Kuranishi family is universal. 
We have analogs to these notions for the transverse holomorphic case, 
we give a brief review of these notions following \cite{GHS}. 

{\it A germ of a deformation of a transverse holomorphic foliation} $\mathcal{F}$ on a manifold $M$ with base space $(B,0)$ is 
given by an open cover $\{U_\alpha \}$ of $M$ and a family of local submersions $f_{\alpha, t}:U_\alpha\rightarrow \C^n$ parametrized 
by $(B,0)$ that are holomorphic in $t\in B$ for each $x\in U_\alpha$. 
If $\Theta_{\mathcal F}$ denotes the sheaf that encodes information about infinitesimal transverse holomorphic deformations
one also has a Kodaira-Spencer map $\rho: T_0B\rightarrow H^1(M,\Theta_{\mathcal F})$ that sends 
$\frac{\partial}{\partial t}$ to certain class in $H^1(M,\Theta_{\mathcal F})$ defined by a section $\theta_{\alpha, \beta}$ 
of the sheaf $\Theta_{\mathcal F}|U_\alpha\cap U_\beta$.  One can consider the full cohomology ring $H^*(M,\Theta_{\mathcal F})$, 
these were proven to be finite dimensional. 

In \cite{GHS} it is showed that there is a versal Kuranishi space of deformations given by the map $\Phi:U\rightarrow H^2(M,\Theta_{\mathcal F})$, 
for $U$ open set in $H^1(M,\Theta_{\mathcal F})$, here, as the complex case, the base of parametrizations is given by $\Phi^{-1}(0)$.
We have that if $H^2(M,\Theta_{\mathcal F})=0$, then the Kuranishi family of deformations of ${\mathcal F}$ 
is isomorphic  to an open set in $H^1(M,\Theta_{\mathcal F})$.  Otherwise, the Kuranishi space is singular. We are not able to prove the vanishing of $H^2(M,\Theta_{\mathcal F})=0$, even so the Kuranishi space is smooth.

\subsection{Deformations of Null Sasaki $\eta$-Einstein in Dimension five} 
From the triviality of the canonical line bundle one obtains the following isomorphisms on a smooth K3 surface $X$: 
$H^0(X, T_X)\cong H^0(X,\Omega_X)=0$, $H^1(X, T_X)\cong H^1(X,\Omega_X)$ and $H^2(X, T_X)\cong H^2(X,\Omega_X)=0.$ 
So the Kuranishi family is universal at all points in a small neighborhood $U$ of the point in the base corresponding to $X$. This base space is smooth of dimension 20. Now, we need to investigate the space of deformations of weighted K3 surfaces $X_{\bf w}$ to find information about the null Sasaki $\eta$-Einstein structures on simply connected five dimensional manifolds. We will do this looking at the minimal resolution $f:\widetilde{X}\rightarrow X_{\bf w}$ of a weighted K3 surface, this resolution contains information that we interpret for 
the null Sasaki $\eta$-Einstein structure in the associated link. 
Recall that in Theorem 3.1,  projectivity of the orbifold is a requirement to obtain orbifold riemannian submersions and therefore Sasakian structures. However it is known that deformations of a projective K3 surface are unstable, that is,  projectivity is not necessarily preserved (unlike the Fano case, in this case we have  $H^2(X, \Omega_X)\not = 0$). 
We will use deformations of polarized K3 surfaces to avoid this lack of stability. 
Before we need the following lemma.

\begin{lemma}
Let $X$ be a Calabi-Yau orbifold with no branch divisor. Denote its minimal resolution by $\widetilde{X}$.  
Then the sublattice $\mathbb Z[\sum_i\Delta_i]\subset H^2(\widetilde{X}, {\mathbb Z})$  generated by the exceptional locus  $\sum_i\Delta_i$ on $\widetilde{X}$ is primitive. 
\end{lemma}
\noindent{\bf Proof.}  Suppose on the contrary that is not primitive, then 
$H^2(\widetilde{X}\backslash\sum_i\Delta_i, \mathbb Z)$ has torsion. 
This implies that $H_1(\widetilde{X}\backslash\sum_i\Delta_i, \mathbb Z)$ has torsion. 
Hence $\pi_1(\widetilde{X}\backslash\sum_i\Delta_i)$ has torsion, and so does  
$\pi_1(X\backslash Sing(X))=\pi_1^{orb}(X)$, obviously a contradiction. \hfill$\blacksquare$
\medskip

\noindent{\bf Theorem 2.} {\it  
Simply connected five dimensional manifolds $N$ that admit null Sasaki $\eta$-Einstein metrics have a space of deformation of dimension equal to $2(b_2(N)-2)$. Moreover, 
the Kuranishi space of this type of  deformations is smooth: 
the moduli space of null Sasaki $\eta$-Einstein structures for $N$ is determined by the following quadric in $\mathbb C\mathbb P^{b_2(N)-1}$: 
$$\{ [\alpha]\in H^2(N, \mathbb C)\, |  \, \rVert([\alpha], [\alpha])\rVert-\rVert([\alpha], [ \overline{\alpha}])\rVert >0\}/\mathbb C^*.$$} 
\noindent{\bf Proof.} First we prove this result at the dimensional level. 
We only need to focus on K3 surfaces with singularities of type $A_{n_i}$ and the corresponding minimal resolution 
$f: X\rightarrow X_{\bf w}$. Here we use the standard procedure:  one notices that at resolving singularities of type $A_{n_i}$ 
the configuration given by the corresponding $(-2)$-curves, denoted by $\sum_i\Delta_i$, creates a sublattice 
$\Z[\sum_i\Delta_i]\subset H^2(X, \Z)$. 
It is known that in any deformation of the minimal resolution $X$ a contraction can be made to yield a deformation on  $X_{\bf w}$. 
So if we consider the local deformations ``normal to the lattice'' $\Z[\sum_i\Delta_i]$ we are obtaining versal deformations of the 
orbifolds $X_{\bf w}$ that keep the orbifold singularities fixed. Additionally, if we consider deformations normal to the lattice spanned 
by $\Z[\sum_i\Delta_i]$ and an ample line 
bundle $L$ one obtains algebraic deformations on the orbifold $X_{\bf w}$ that fix the singularities 
(in analogy to the polarization of non-singular K3 surfaces). This type of deformation is the one we are interested in since this is equivalent 
to deformations keeping the weighted $S^1$ action on the Seifert bundle $S^1\hookrightarrow N \rightarrow X_{\bf w}$ fixed. 

Let us denote by $M$ the lattice spanned  by $L$ and $\Z[\sum_i\Delta_i]$. From Lemma 4.1, $M$ is a primitive sublattice of the 
K3 lattice $L_{\Lambda}$. It has signature $(1, \mathrm{rank} M-1)$. Let $T=M^\perp$, one can associate to it the period domain 
of marked $M$-polarized K3 surfaces 
\begin{equation}
\Omega_M=\{[\omega]\in {\mathbb P}(T\otimes \mathbb C)|\, \omega.\omega=0,\, \omega.\bar{\omega}>0  \}.
\end{equation}
By the Torrelli theorem for K3 surfaces, $\Omega_M$ is the moduli space of marked K3 surfaces $(X,\phi)$, with 
$\phi:H^2(X,\Z)\rightarrow L_{\Lambda}$ a marking such that $\phi^{-1}(M)\subseteq \mathrm{NS} (X)$. 
In \cite{Dol} Dolgachev  proves that the quotient of $\Omega_M$ via the group action of elements $\varphi\in \mathrm{Aut} (L_{\Lambda})$ 
such that $\varphi(M)=M$ is the moduli space of $M$-polarized K3 surfaces. In particular, 
the dimension of this moduli space is $20-$rank$(M)$. 

Returning to our problem, now it is straightforward to give a formula for the dimensions of our moduli spaces in terms of their 
second Betti numbers $b_2(X_{\bf w})$ since the rank of the lattice $M$ is directly related with this number. For Seifert bundles $S^1\hookrightarrow N\rightarrow X$, Boyer and Galicki showed (see \cite{BBG} Proposition 8.2.6)  the existence of 
an exact sequence 
$$0\rightarrow H^1(X_{\bf w}, T_{X_{\bf w}})\rightarrow H^1(N,\Theta_{\mathcal F})\rightarrow 
H^0(X_{\bf w}, T_{X_{\bf w}})\rightarrow H^2(X_{\bf w}, T_{X_{\bf w}}).$$
Hence, under the absence of infinitesimal holomorphic automorphism on the orbifold, all the deformations of the transverse holomorphic structure come from deforming the complex structure on $X$. Since on a K3 surface  
$H^0(X_{\bf w},T_{X_{\bf w}})=0$ one has $$H^1(N,\Theta_{\mathcal F})=H^1(X_{\bf w},T_{X_{\bf w}}).$$ 
Combined with Theorem 3.2,  one obtains that the (real) dimension of the space of deformations on the corresponding total 
space via the  Seifert bundle has dimension $2(b_2(N)-2)$.  

But we can be more precise. 
As stated before, on a compact Sasakian manifold (although K-contact suffices) there is a transverse Hodge theory. The transverse Hodge star 
is defined in terms of the usual Hodge star by 
$$\bar{\star}\alpha=\star(\eta\wedge \alpha)=(-1)^k\xi\rfloor \star \alpha.$$ There is also a {\it basic Laplacian} $\Delta_B$ 
defined in terms of $d_B$ and its adjoint. The kernel of $\Delta_B$ is the space 
$\mathcal{H}^k_B(\mathcal{F}_\xi)$ of basic harmonic $k$ forms. The transverse Hodge theorem \cite{EKA} gives the isomorphism 
$\mathcal{H}^k_B(\mathcal{F}_\xi)=H^k_B(\mathcal{F}_\xi)$. Since $\mathcal{F}_\xi$ is {\it taut} (for example, when $\eta$ is invariant under the flow of $\xi$, a condition that is mandatory on compact Sasakian manifolds),  
the pairing $$\Omega_B^k\otimes \Omega_B^{2n-k}\longrightarrow \mathbb R$$ given by 
$\alpha\otimes \beta \mapsto \int_N\eta\wedge (\alpha\wedge \beta)$ induces an isomorphism on the transverse harmonic forms and it defines 
a non-degenerate inner product $\langle \alpha, \beta \rangle_B=\int_N\eta \wedge\alpha\wedge \bar{\star}\beta$ (see \cite{Ton} for details). 
Hence  Poincar\'e duality in $H^*_B(\mathcal{F}_\xi)$ holds, 
that is, it induces a non-degenerate pairing $H^k_B(\mathcal{F}_\xi)\times H_B^{2n-k}(\mathcal{F}_\xi)\rightarrow \mathbb R$ that 
from now on we will denote also by $( [\alpha]_B, [\beta]_B).$ 
It is known (see \cite{BBG} equation 7.2.1) that there is the Gysin sequence that relates basic and De Rham comology, combining this with the fact 
that here we are dealing with orbifolds such that $H^1_B(\mathcal{F}_\xi)=0$ one has 
$$
0\longrightarrow H_B^0(\mathcal{F}_\xi)\stackrel{\delta}\rightarrow H_B^2(\mathcal{F}_\xi) \longrightarrow H^2(N, \mathbb R)\longrightarrow 0,  
$$
where $\delta$ is the conecting homomorphism given by $\delta[\alpha]=[d\eta\wedge \alpha]_B$. 
Hence $H_B^2(\mathcal{F}_\xi)=H_B^0(\mathcal{F}_\xi)\oplus H^2(N, \mathbb R)$. Thus, one can write any class $[\alpha]_B\in H_B^2(\mathcal{F}_\xi)$ as $[\alpha]_B=z[d\eta]_B +[\alpha]$ where 
$[\alpha]\in H^2(N, \mathbb R)$ and $z\in \mathbb C$. 

We extend this result by linearity over the complex numbers and one obtains 
$([\alpha]_B, [\alpha]_B)=z_1z_2([d\eta]_B, [d\eta]_B)+([\alpha],[\beta])$ with $z_1, z_2\in \mathbb C$. 
The following two identities follow immediately.
\begin{eqnarray}
 ([\alpha]_B,[\alpha]_B) &=& ([\alpha],[\alpha])+ z^2([d\eta]_B,[d\eta]_B)\\
([\alpha]_B,[\overline\alpha]_B) &=& ([\alpha],[\overline\alpha])+ z\overline{z}([d\eta]_B,[d\eta]_B). 
\end{eqnarray}
Here $([d\eta]_B,[d\eta]_B) \not=0$ since $[d\eta]_B\in H_B^2(\mathcal{F}_\xi)$ is non-trivial. 
Thus we obtain a non-degenerate quadratic form on $H^2(N, \mathbb C)$. 
Since $H^1(N,\Theta_{\mathcal F})=H^1(X_{\bf w},T_{X_{\bf w}})$, to obtain information about the moduli space on null Sasaki $\eta$-Einstein structures, it is enough to look at the information about the period domain described in $(4)$: what happens to elements in $H^2(N, \mathbb C)$ if {\bf a)} $([\alpha]_B,[\alpha]_B)=0$ 
and {\bf b)} $([\alpha]_B,[\overline\alpha]_B)>0$? 
Equations $(5)$ and $(6)$ give us the possibility to express this conditions in terms of the corresponding classes in $H^2(N, \mathbb C)$. 
Elementary calculations yield to the open set given in the statement (notice that since we are considering $[\alpha]_B\in H^2_B(\mathcal{F}_\xi)$ an element in the period space, it is non-trivial, 
which implies that its corresponding 
representative in $[\alpha]\in H^2(N, \mathbb C)$ is non-trivial as well).  
\hfill$\blacksquare$
\medskip

The crucial ideas used in the proof of Theorem 2 are very special to the case when the base space is a K3 surface, we expand on the issues in the following two remarks. 

\begin{remark}{\rm 
In the non-simply connected case one has the same isomorphism $H^1(N,\Theta_{\mathcal F})=H^1(X,T_{X}),$  when $N$ 
is an abelian variety (where it is known there are not infinitesimal automorphism). So the Kuranishi family is non-singular. 
However, this time we cannot find expressions similar to the ones given in equations $(4)$ and $(5)$ (the lack of simply connectivity does not allow it)}. 
\end{remark}

\begin{remark}{\rm  Let $S^1\hookrightarrow N\rightarrow Z$ be a Seifert bundle where $S^1$  acts locally free at worst, $N$ is a 
$(4n+1)$-dimensional manifold and $Z$ is a simply connected hyperk\"ahler orbifold. 
It is known \cite{Huy}, as in the case of K3 surfaces, that there is a non-degenerate quadratic form in $H^2(Z, \mathbb Z)$ when $Z$ is a manifold. 
The period map is surjective and locally an isomorphism. 
By $(4)$ we can use this form to obtain a non-degenerate quadratic form on $H^2(N,\mathbb C)$ and we have an identical result 
as the one stated in the previous corollary provided the existence of a Global Torelli theorem, this is not always the case, see \cite{Huy}.  
In the orbifold situation, 
Tondeur's result \cite{Ton} on non-degeneracy of the quadratic form  extends 
this form on orbifold hyperk\"ahler structures. However, here the idea of polarized lattices that, in the K3 orbifold case, 
allowed having period maps with the properties stated above, is not available in general (moreover, the existence of crepant resolutions 
is not even known 
in higher dimensions).}   
\end{remark}

\bigskip

In \cite{BBG}, Theorem 5.5.7  a formula is given to calculate the real dimension of (polarized) K\"ahler-Einstein metrics 
up to homothety in each positive (1,1) cohomology class for quasi-smooth  weighted hypersurfaces $X_f$ of degree $d$ 
in ${\mathbb P}({\bf w})$ with $|{\bf w}|-d=0$: 
$$2[h^0({\mathbb P}({\bf w}), {\mathcal{O}}(d))-\sum_i h^0({\mathbb P}({\bf w}),{\mathcal{O}}(w_i))],$$ here $h^0({\mathbb P}({\bf w}), {\mathcal{O}}(l))$ denotes, as usual, the dimension of the space of weighted homogeneous 
polynomials of degree $l$. 
The theorem proved above, gives a different and more general way (not only for links of hypersurfaces) to find the dimension of the moduli of simply connected 
null Sasaki $\eta$-Einstein structures in terms of the second Betti number. For example, consider number 2 in Table 1: $X_5\in \mathbb P(1,1,1,2)$ which corresponds, via Seifert bundle, to $\# 20 (S^2\times S^3)$. The general polynomial $f({\bf z})$ defining this hypersurface can be written as the sum 
$f({\bf z})=g^{(5)}(z_0, z_1, z_2)+g^{(3)}(z_0, z_1, z_2)z_3 + g^{(1)}(z_0, z_1, z_2)z_3^2,$ where $g^{(i)}(z_0, z_1, z_2)$ are 
homogeneous polynomials of degree $i$. 
Hence, $\displaystyle{h^0({\mathbb P}({\bf w}), {\mathcal{O}}(d))= {7\choose 2}+ {5\choose 2}+{3\choose 2}=34.}$ 
On the other hand, $H^0({\mathbb P}({\bf w}), {\mathcal{O}}(w_i))$ is defined on generators  by 
$$\left ( \begin{array}{cc}
{\mathbb A} \left (\begin{array}{c}
    z_0\\
           z_1\\
           z_2
                  
\end{array}
\right )\\

\alpha z_3 +\varphi^{(2)}(z_0, z_1, z_2)
 
\end{array}
   \right ),$$ 
where $\mathbb A\in GL(3, \mathbb C)$, $\alpha\in \C^*$ and $\varphi^{(2)}$ is any homogeneous polynomial of degree 2 in $(z_0, z_1, z_2)$. 
Thus $\sum_i h^0({\mathbb P}({\bf w}), {\mathcal{O}}(w_i))= 9+1+{4\choose 2}=16.$ So we obtain that the dimension of the moduli of (polarized) K\"ahler-Einstein metrics (and therefore the moduli of 
null Sasaki $\eta$-Einstein structures) is 36 which coincides with $2(k-2)$, for $\# 20 (S^2\times S^3).$ Similar calculations yield to 
the same dimension for the K3 weighted hypersurface of degree 12 in $\mathbb P(1,1,4,6)$  (number 14 in Table 1). 
This was expected since $X_{12}\in \mathbb P(1,1,4,6)$ and $X_5\in \mathbb P(1,1,1,2)$ both have second Betti number 21. 
Moreover, number 2 in Table 2 $X_{3,3}\subset {\mathbb P}(1,1,1,1,2)$ has the same moduli number since its second Betti number is 21 as well.




\section{Tables}

\tiny{
\noindent{\bf Table 1. Reid's List of 95 Codimension 1 Weighted K3 Surfaces.}
$$
\begin{array}{|c|c|c|c|}
\hline
    & X_{\bf w}                & {\rm singularities}  & b_2(X_{\bf w}) \\ \hline 
\mathrm{No.}1 & X_4\subset \mathbb P (1,1,1,1)  &      &  22   \\  \hline
\mathrm{No.}2 & X_5\subset \mathbb P (1,1,1,2)  &  A_1 &  21  \\  \hline
\mathrm{No.}3 & X_6\subset \mathbb P (1,1,1,3)  &      &  22    \\  \hline
\mathrm{No.}4 & X_6\subset \mathbb P (1,1,2,2)  &  3\times A_1 & 19\\  \hline
\mathrm{No.}5 & X_7\subset \mathbb P (1,1,2,3)  &  A_1,A_2     & 19    \\  \hline
\mathrm{No.}6 & X_8\subset \mathbb P (1,1,2,4)  &   2\times A_1 & 20    \\  \hline
\mathrm{No.}7 & X_8\subset \mathbb P (1,2,2,3)  &   4\times A_1, A_2 & 16\\  \hline
\mathrm{No.}8 & X_9\subset \mathbb P (1,1,3,4)  &    A_3 & 19 \\  \hline
\mathrm{No.}9 & X_9\subset \mathbb P (1,2,3,3)  &   A_1, 3\times A_2 & 15\\  \hline
\mathrm{No.}10 & X_{10}\subset \mathbb P (1,1,3,5)  & A_2  & 20 \\  \hline
\mathrm{No.}11 & X_{10}\subset \mathbb P (1,2,2,5)  &  5\times A_1 & 17\\  \hline
\mathrm{No.}12 & X_{10}\subset \mathbb P (1,2,3,4)  &  2\times A_1, A_2, A_3 & 15\\  \hline
\mathrm{No.}13 & X_{11}\subset \mathbb P (1,2,3,5)  &   A_1, A_2, A_5 & 15 \\  \hline
\mathrm{No.}14 & X_{12}\subset \mathbb P (1,1,4,6)  &  A_1  & 21\\  \hline
\mathrm{No.}15 & X_{12}\subset \mathbb P (1,2,3,6)  &  2\times A_1, 2\times A_2 & 14\\  \hline
\mathrm{No.}16 & X_{12}\subset \mathbb P (1,2,4,5)  &   3\times A_1, A_4 & 15\\  \hline
\mathrm{No.}17 & X_{12}\subset \mathbb P (1,3,4,4)  &   3\times A_3 & 13\\  \hline
\mathrm{No.}18 & X_{12}\subset \mathbb P (2,2,3,5)  &   6\times A_1, A_4  & 12\\  \hline
\mathrm{No.}19 & X_{12}\subset \mathbb P (2,3,3,4)  &   3\times A_1, 4\times A_2 & 11\\  \hline
\mathrm{No.}20 & X_{13}\subset \mathbb P (1,3,4,5)  &   A_2, A_3, A_4 & 13\\  \hline
\mathrm{No.}21 & X_{14}\subset \mathbb P (1,2,4,7)  &   3\times A_1, A_3  & 16 \\  \hline
\mathrm{No.}22 & X_{14}\subset \mathbb P (2,2,3,7)  &   7\times A_1, A_2 & 13\\  \hline
\mathrm{No.}23 & X_{14}\subset \mathbb P (2,3,4,5)  &   3\times A_1, A_2, A_3, A_4 & 10\\  \hline
\mathrm{No.}24 & X_{15}\subset \mathbb P (1,2,5,7)  &   A_1, A_6 & 15 \\  \hline
\mathrm{No.}25 & X_{15}\subset \mathbb P (1,3,4,7)  &    A_3, A_6 & 13\\  \hline
\mathrm{No.}26 & X_{15}\subset \mathbb P (1,3,5,6)  &   2\times A_2, A_5 & 13\\  \hline
\mathrm{No.}27 & X_{15}\subset \mathbb P (2,3,5,5)  &   A_1, 3\times A_4 & 9\\  \hline
\mathrm{No.}28 & X_{15}\subset \mathbb P (3,3,4,5)  &   5\times A_2, A_3 & 9 \\  \hline
\mathrm{No.}29 & X_{16}\subset \mathbb P (1,2,5,8)  &   2\times A_1, A_4 & 16 \\  \hline
\mathrm{No.}30 & X_{16}\subset \mathbb P (1,3,4,8)  &    A_2, 2\times A_3 & 14\\  \hline
\mathrm{No.}31 & X_{16}\subset \mathbb P (1,4,5,6)  &    A_1, A_4, A_5   & 12\\  \hline
\mathrm{No.}32 & X_{16}\subset \mathbb P (2,3,4,7)  &   4\times A_1, A_2, A_6 & 10\\  \hline
\mathrm{No.}33 & X_{17}\subset \mathbb P (2,3,5,7)  &    A_1, A_2, A_4, A_6 & 9\\  \hline
\mathrm{No.}34 & X_{18}\subset \mathbb P (1,2,6,9)  &   3\times A_1, A_2 & 15\\  \hline
\mathrm{No.}35 & X_{18}\subset \mathbb P (1,3,5,9)  &   2\times A_2, A_4 & 14 \\  \hline
\mathrm{No.}36 & X_{18}\subset \mathbb P (1,4,6,7)  &  A_3, A_1, A_6  & 12 \\  \hline
\mathrm{No.}37 & X_{18}\subset \mathbb P (2,3,4,9)  &  4\times A_1, 2\times A_2, A_3 & 11\\  \hline
\mathrm{No.}38 & X_{18}\subset \mathbb P (2,3,5,8)  &  2\times A_1, A_4, A_7 & 9 \\  \hline
\mathrm{No.}39 & X_{18}\subset \mathbb P (3,4,5,6)  &  3\times A_2, A_3, A_1, A_4 & 8 \\  \hline
\mathrm{No.}40 & X_{19}\subset \mathbb P (3,4,5,7)  &  A_2, A_3, A_4, A_6 & 7 \\  \hline
\mathrm{No.}41 & X_{20}\subset \mathbb P (1,4,5,10)  &  A_1, 2\times A_4 & 13 \\  \hline
\mathrm{No.}42 & X_{20}\subset \mathbb P (2,3,5,10)  &  2\times A_1, A_2, 2\times A_4 & 10 \\  \hline
\mathrm{No.}43 & X_{20}\subset \mathbb P (2,4,5,9)  &  5\times A_1, A_8  & 9 \\  \hline
\mathrm{No.}44 & X_{20}\subset \mathbb P (2,5,6,7)  &  3\times A_1, A_5, A_6 & 8 \\  \hline
\mathrm{No.}45 & X_{20}\subset \mathbb P (3,4,5,8)  &  A_2, 2\times A_3, A_7 & 7 \\  \hline
\mathrm{No.}46 & X_{21}\subset \mathbb P (1,3,7,10)  &  A_9 & 13 \\  \hline
\mathrm{No.}47 & X_{21}\subset \mathbb P (1,5,7,8)  &  A_4, A_7 & 11\\  \hline
\mathrm{No.}48 & X_{21}\subset \mathbb P (2,3,7,9)  &  A_1, 2\times A_2, A_8 & 9 \\  \hline
\mathrm{No.}49 & X_{21}\subset \mathbb P (3,5,6,7)  &  3\times A_2, A_4, A_5  & 7 \\  \hline
\mathrm{No.}50 & X_{22}\subset \mathbb P (1,3,7,11)  &  A_2, A_6  & 14 \\  \hline
\mathrm{No.}51 & X_{22}\subset \mathbb P (1,4,6,11)  &  A_3, A_1, A_5 & 13\\  \hline
\mathrm{No.}52 & X_{22}\subset \mathbb P (2,4,5,11)  &  5\times A_1, A_3, A_4 & 10\\  \hline
\mathrm{No.}53 & X_{24}\subset \mathbb P (1,3,8,12)  &  2\times A_2, A_3 & 16\\  \hline
\mathrm{No.}54 & X_{24}\subset \mathbb P (1,6,8,9)  &  A_1, A_2, A_8 & 11\\  \hline
\mathrm{No.}55 & X_{24}\subset \mathbb P (2,3,7,12)  &  2\times A_1, 2\times A_2, A_6 & 10\\  \hline
\mathrm{No.}56 & X_{24}\subset \mathbb P (2,3,8,11)  &  3\times A_1, A_{10} & 9\\  \hline
\mathrm{No.}57 & X_{24}\subset \mathbb P (3,4,5,12)  &  2\times A_2, 2\times A_3, A_4 & 8\\  \hline
\mathrm{No.}58 & X_{24}\subset \mathbb P (3,4,7,10)  &  A_1, A_6, A_9 &  6 \\  \hline
\mathrm{No.}59 & X_{24}\subset \mathbb P (3,6,7,8)  &  4\times A_2, A_1, A_6 & 7 \\  \hline
\mathrm{No.}60 & X_{24}\subset \mathbb P (4,5,6,9)  &  2\times A_1, A_4, A_2, A_8 & 6 \\  \hline
\end{array}
$$

\noindent{\bf Table 1. Reid's List of 95 Codimension 1 Weighted K3 Surfaces.}
$$
\begin{array}{|c|c|c|c|}
\hline
    & X_{\bf w}                & {\rm singularities}  & b_2(X_{\bf w}) \\ \hline 
\mathrm{No.}61 & X_{25}\subset \mathbb P (4,5,7,9)  &  A_3, A_6, A_8 & 5 \\  \hline
\mathrm{No.}62 & X_{26}\subset \mathbb P (1,5,7,13)  &  A_4, A_6 & 12 \\  \hline
\mathrm{No.}63 & X_{26}\subset \mathbb P (2,3,8,13)  &  3\times A_1, A_2, A_7 & 10\\  \hline
\mathrm{No.}64 & X_{26}\subset \mathbb P (2,5,6,13)  &  4\times A_1, A_4, A_5 & 9 \\  \hline
\mathrm{No.}65 & X_{27}\subset \mathbb P (2,5,9,11)  &  A_1, A_4, A_{10} & 7 \\  \hline
\mathrm{No.}66 & X_{27}\subset \mathbb P (5,6,7,8)  &  A_4, A_5, A_2, A_6 & 5 \\  \hline
\mathrm{No.}67 & X_{28}\subset \mathbb P (1,4,9,14)  &  A_1, A_8 & 13 \\  \hline
\mathrm{No.}68 & X_{28}\subset \mathbb P (3,4,7,14)  &  A_2, A_1, 2\times A_6 & 7 \\  \hline
\mathrm{No.}69 & X_{28}\subset \mathbb P (4,6,7,11)  &  2\times A_1, A_5, A_{10} & 5 \\  \hline
\mathrm{No.}70 & X_{30}\subset \mathbb P (1,4,10,15)  &  A_3, A_4, A_1 & 14 \\  \hline
\mathrm{No.}71 & X_{30}\subset \mathbb P (1,6,8,15)  &  A_1, A_2, A_7 & 12 \\  \hline
\mathrm{No.}72 & X_{30}\subset \mathbb P (2,3,10,15)  &  3\times A_1, 2\times A_2, A_4 & 6 \\  \hline
\mathrm{No.}73 & X_{30}\subset \mathbb P (2,6,7,15)  &  5\times A_1, A_2, A_6 & 9 \\  \hline
\mathrm{No.}74 & X_{30}\subset \mathbb P (3,4,10,13)  &  A_3, A_1, A_{12} & 5 \\  \hline
\mathrm{No.}75 & X_{30}\subset \mathbb P (4,5,6,15)  &  A_3, 2\times A_1, 2\times A_4, A_2 & 7\\  \hline
\mathrm{No.}76 & X_{30}\subset \mathbb P (5,6,8,11)  &  A_1, A_7, A_{10} & 4\\  \hline
\mathrm{No.}77 & X_{32}\subset \mathbb P (2,5,9,16)  &  2\times A_1, A_4, A_8 & 8\\  \hline
\mathrm{No.}78 & X_{32}\subset \mathbb P (4,5,7,16)  &  2\times A_3, A_4, A_6 & 6\\  \hline
\mathrm{No.}79 & X_{33}\subset \mathbb P (3,5,11,14)  &  A_4, A_{13} & 5\\  \hline
\mathrm{No.}80 & X_{34}\subset \mathbb P (3,4,10,17)  &  A_2, A_3, A_1, A_9 & 7\\  \hline
\mathrm{No.}81 & X_{34}\subset \mathbb P (4,6,7,17)  &  A_3, 2\times A_1, A_5, A_6 & 6\\  \hline
\mathrm{No.}82 & X_{36}\subset \mathbb P (1,5,12,18)  &  A_4, A_5 & 13 \\  \hline
\mathrm{No.}83 & X_{36}\subset \mathbb P (3,4,11,18)  &  2\times A_2, A_1, A_{10}& 7 \\  \hline
\mathrm{No.}84 & X_{36}\subset \mathbb P (7,8,9,12)  &  A_6, A_7, A_3, A_2 & 4 \\  \hline
\mathrm{No.}85 & X_{38}\subset \mathbb P (3,5,11,19)  &  A_2, A_4, A_{10} & 6 \\  \hline
\mathrm{No.}86 & X_{38}\subset \mathbb P (5,6,8,19)  &  A_4, A_5, A_1, A_7 & 5\\  \hline
\mathrm{No.}87 & X_{40}\subset \mathbb P (5,7,8,20)  &  2\times A_4, A_6, A_3 & 5 \\  \hline
\mathrm{No.}88 & X_{42}\subset \mathbb P (1,6,14,21)  &  A_1, A_2, A_6 & 13 \\  \hline
\mathrm{No.}89 & X_{42}\subset \mathbb P (2,5,14,21)  &  3\times A_1, A_4, A_6 & 9 \\  \hline
\mathrm{No.}90 & X_{42}\subset \mathbb P (3,4,14,21)  &  2\times A_2, A_3, A_1, A_6 & 8 \\  \hline
\mathrm{No.}91 & X_{44}\subset \mathbb P (4,5,13,22)  &  A_1, A_4, A_{12} & 5 \\  \hline
\mathrm{No.}92 & X_{48}\subset \mathbb P (3,5,16,24)  &  2\times A_2, A_4, A_7 & 7 \\  \hline
\mathrm{No.}93 & X_{50}\subset \mathbb P (7,8,10,25)  &  A_6, A_7, A_1, A_4 & 4 \\  \hline
\mathrm{No.}94 & X_{54}\subset \mathbb P (4,5,18,27)  &  A_3, A_1, A_4, A_8 & 6 \\  \hline
\mathrm{No.}95 & X_{66}\subset \mathbb P (5,6,22,33)  &  A_4, A_1, A_2, A_{10} & 5\\  \hline
\end{array}
$$


\medskip

\noindent{\bf Tabel 2. Fletcher's List of 84 Codimension 2 Weighted K3 surfaces.}
$$
\begin{array}{|c|c|c|c|}
\hline
    & X_{\bf w}                & {\rm singularities}       & b_2(X_{\bf w})\\ \hline 

\mathrm{No.}1 & X_{2,3}\subset \mathbb P (1,1,1,1,1)  &   & 22   \\  \hline
\mathrm{No.}2 & X_{3,3}\subset \mathbb P (1,1,1,1,2)  &   A_1 & 21  \\  \hline
\mathrm{No.}3 & X_{3,4}\subset \mathbb P (1,1,1,2,2)  &   2\times A_1 & 20 \\  \hline
\mathrm{No.}4 & X_{4,4}\subset \mathbb P (1,1,1,2,3)  &  A_2  & 20  \\  \hline
\mathrm{No.}5 & X_{4,4}\subset \mathbb P (1,1,2,2,2)  &  4\times A_1   & 18 \\  \hline
\mathrm{No.}6 & X_{4,5}\subset \mathbb P (1,1,2,2,3)  &   2\times A_1, A_2 & 18 \\  \hline
\mathrm{No.}7 & X_{4,6}\subset \mathbb P (1,1,2,3,3)  &   2\times A_2   & 18 \\  \hline
\mathrm{No.}8 & X_{4,6}\subset \mathbb P (1,2,2,2,3)  &    6\times A_1   & 16\\  \hline
\mathrm{No.}9 & X_{5,6}\subset \mathbb P (1,1,2,3,4)  &   A_1, A_3  & 18 \\   \hline
\mathrm{No.}10 & X_{5,6}\subset \mathbb P (1,2,2,3,3)  & 3\times A_1, 2\times A_2 & 15 \\  \hline
\mathrm{No.}11 & X_{6,6}\subset \mathbb P (1,1,2,3,5)  &  A_4  & 18 \\  \hline
\mathrm{No.}12 & X_{6,6}\subset \mathbb P (1,2,2,3,4)  &  4\times A_1, A_3 & 15 \\  \hline
\mathrm{No.}13 & X_{6,6}\subset \mathbb P (1,2,3,3,3)  &   4\times A_2 & 14 \\  \hline
\mathrm{No.}14 & X_{6,6}\subset \mathbb P (2,2,2,3,3)  &  9\times A_1 & 13\\  \hline
\mathrm{No.}15 & X_{6,7}\subset \mathbb P (1,2,2,3,5)  &  3\times A_1, A_4  & 15 \\  \hline
\mathrm{No.}16 & X_{6,7}\subset \mathbb P (1,2,3,3,4)  &   A_1, 2\times A_2, A_3 & 14 \\  \hline
\mathrm{No.}17 & X_{6,8}\subset \mathbb P (1,1,3,4,5)  &   A_4  & 18  \\  \hline
\mathrm{No.}18 & X_{6,8}\subset \mathbb P (1,2,2,3,5)  &   2\times A_2, A_4  & 14\\  \hline
\mathrm{No.}19 & X_{6,8}\subset \mathbb P (1,2,3,4,4)  &   2\times A_1, 2\times A_3 & 14\\  \hline
\mathrm{No.}20 & X_{6,8}\subset \mathbb P (2,2,3,3,4)  &   6\times A_1, 2\times A_2 & 12\\  \hline
\mathrm{No.}21 & X_{6,9}\subset \mathbb P (1,2,3,4,5)  &   A_1, A_3, A_4   & 14 \\  \hline
\mathrm{No.}22 & X_{7,8}\subset \mathbb P (1,2,3,4,5)  &   2\times A_1, A_2, A_4 & 14\\  \hline
\mathrm{No.}23 & X_{6,10}\subset \mathbb P (1,2,3,5,5)  &   2\times A_4 & 14 \\  \hline
\mathrm{No.}24 & X_{6,10}\subset \mathbb P (2,2,3,4,5)  &   7\times A_1, A_3 & 12\\  \hline
\mathrm{No.}25 & X_{8,9}\subset \mathbb P (1,2,3,4,7)  &    2\times A_1, A_6 & 14 \\  \hline
\mathrm{No.}26 & X_{8,9}\subset \mathbb P (1,3,4,4,5)  &   2\times A_3, A_4 & 13\\  \hline
\mathrm{No.}27 & X_{8,9}\subset \mathbb P (2,3,3,4,5)  &  2\times A_1, 3\times A_2, A_4 & 10 \\  \hline
\mathrm{No.}28 & X_{8,10}\subset \mathbb P (1,2,3,5,7)  &   A_2, A_6 & 14\\  \hline
\mathrm{No.}29 & X_{8,10}\subset \mathbb P (1,2,4,5,6)  &   3\times A_1, A_5 & 14\\  \hline
\mathrm{No.}30 & X_{8,10}\subset \mathbb P (1,3,4,5,5)  &    A_2, 2\times A_4 & 12 \\  \hline
\mathrm{No.}31 & X_{8,10}\subset \mathbb P (2,3,4,4,5)  &    4\times A_1, A_2, 2\times A_3  & 10\\  \hline
\mathrm{No.}32 & X_{9,10}\subset \mathbb P (1,2,3,5,8)  &   A_1, A_7  & 14\\  \hline
\mathrm{No.}33 & X_{9,10}\subset \mathbb P (1,3,4,5,6)  &    A_2, A_3, A_5 & 12 \\  \hline
\mathrm{No.}34 & X_{9,10}\subset \mathbb P (2,2,3,5,7)  &   5\times A_1, A_6 & 11\\  \hline
\mathrm{No.}35 & X_{9,10}\subset \mathbb P (2,3,4,5,5)  &   2\times A_1, A_3, 2\times A_4 & 9  \\  \hline
\end{array}
$$

\noindent{\bf Tabel 2. Fletcher's List of 84 Codimension 2 Weighted K3 surfaces.}
$$
\begin{array}{|c|c|c|c|}
\hline
    & X_{\bf w}                & {\rm singularities}       & b_2(X_{\bf w})\\ \hline 
\mathrm{No.}36 & X_{8,12}\subset \mathbb P (1,3,4,5,7)  &  A_4, A_6 & 12 \\  \hline
\mathrm{No.}37 & X_{8,12}\subset \mathbb P (2,3,4,5,6)  &  4\times A_1, 2\times A_2, A_4 & 10\\  \hline
\mathrm{No.}38 & X_{9,12}\subset \mathbb P (2,3,4,5,7)  &  3\times A_1, A_4, A_6 & 9 \\  \hline
\mathrm{No.}39 & X_{10,11}\subset \mathbb P (2,3,4,5,7)  &  2\times A_1, A_2, A_3, A_6 & 9\\  \hline
\mathrm{No.}40 & X_{10,12}\subset \mathbb P (1,3,4,5,9)  &  A_2, A_8 & 12 \\  \hline
\mathrm{No.}41 & X_{10,12}\subset \mathbb P (1,3,5,6,7)  &  2\times A_2, A_6 & 12 \\  \hline
\mathrm{No.}42 & X_{10,12}\subset \mathbb P (1,2,5,6,6)  &  A_1, 2\times A_5 & 11\\  \hline
\mathrm{No.}43 & X_{10,12}\subset \mathbb P (2,3,4,5,8)  &  3\times A_1, A_3, A_7 & 9 \\  \hline
\mathrm{No.}44 & X_{10,12}\subset \mathbb P (2,3,5,5,7)  &  2\times A_4, A_6 & 8 \\  \hline
\mathrm{No.}45 & X_{10,12}\subset \mathbb P (2,4,5,5,6)  &  5\times A_1, 2\times A_4 & 9 \\  \hline
\mathrm{No.}46 & X_{10,12}\subset \mathbb P (3,3,4,5,7)  &  4\times A_2, A_6 & 8 \\  \hline
\mathrm{No.}47 & X_{10,12}\subset \mathbb P (3,4,4,5,6)  &  2\times A_2, 3\times A_3, A_1 & 8\\  \hline
\mathrm{No.}48 & X_{11,12}\subset \mathbb P (1,4,5,6,7)  &  A_1, A_4, A_6 & 11\\  \hline
\mathrm{No.}49 & X_{10,14}\subset \mathbb P (1,2,5,7,9)  &  A_8   & 14 \\  \hline
\mathrm{No.}50 & X_{10,14}\subset \mathbb P (2,3,5,7,7)  &  A_2, 2\times A_6 & 8\\  \hline
\mathrm{No.}51 & X_{10,14}\subset \mathbb P (2,4,5,6,7)  &  5\times A_1, A_3, A_5 & 9\\  \hline
\mathrm{No.}52 & X_{10,15}\subset \mathbb P (2,3,5,7,8)  &  A_1, A_6, A_7 & 8 \\  \hline
\mathrm{No.}53 & X_{12,13}\subset \mathbb P (3,4,5,6,7)  &  2\times A_2, A_1,A_4, A_6 & 7\\  \hline
\mathrm{No.}54 & X_{12,14}\subset \mathbb P (1,3,4,7,11)  &  A_{10}  & 12 \\  \hline
\mathrm{No.}55 & X_{12,14}\subset \mathbb P (1,4,6,7,8)  &  A_1, A_3, A_7 & 11\\  \hline
\mathrm{No.}56 & X_{12,14}\subset \mathbb P (2,3,4,7,10)  &  4\times A_1, A_9 & 9\\  \hline
\mathrm{No.}57 & X_{12,14}\subset \mathbb P (2,3,5,7,9)  &  A_2, A_4, A_8 & 8\\  \hline
\mathrm{No.}58 & X_{12,14}\subset \mathbb P (3,4,5,7,7)  &  A_4, 2\times A_6 & 6 \\  \hline
\mathrm{No.}59 & X_{12,14}\subset \mathbb P (4,4,5,6,7)  &  3\times A_3, 2\times A_1, A_4 & 7 \\  \hline
\mathrm{No.}60 & X_{12,15}\subset \mathbb P (1,4,5,6,11)  &  A_1, A_{10} & 11\\  \hline
\mathrm{No.}61 & X_{12,15}\subset \mathbb P (3,4,5,6,9)  &  3\times A_2, A_1, A_8 & 7 \\  \hline
\mathrm{No.}62 & X_{12,15}\subset \mathbb P (3,4,5,7,8)  &  A_3, A_6, A_7 & 6 \\  \hline
\mathrm{No.}63 & X_{12,16}\subset \mathbb P (2,5,6,7,8)  &  4\times A_1, A_4, A_6 & 8 \\  \hline
\mathrm{No.}64 & X_{14,15}\subset \mathbb P (2,3,5,7,12)  &  A_1, A_2, A_{11} & 8\\  \hline
\mathrm{No.}65 & X_{14,15}\subset \mathbb P (2,5,6,7,9)  &  2\times A_1, A_5, A_8 & 7  \\  \hline
\mathrm{No.}66 & X_{14,15}\subset \mathbb P (3,4,5,7,10)  &  A_3, A_4, A_9 & 6 \\  \hline
\mathrm{No.}67 & X_{14,15}\subset \mathbb P (3,5,6,7,8)  &  2\times A_2, A_5, A_7 & 6 \\  \hline
\mathrm{No.}68 & X_{14,16}\subset \mathbb P (1,5,7,8,9)  &  A_4, A_8 & 10 \\  \hline
\mathrm{No.}69 & X_{14,16}\subset \mathbb P (3,4,5,7,11)  &  A_2, A_4, A_{10} & 6 \\  \hline
\mathrm{No.}70 & X_{14,16}\subset \mathbb P (4,5,6,7,8)  &  A_1, 2\times A_3, A_4, A_5 & 6 \\  \hline
\mathrm{No.}71 & X_{15,16}\subset \mathbb P (2,3,5,8,13)  &  2\times A_1, A_{12} & 8 \\  \hline
\mathrm{No.}72 & X_{15,16}\subset \mathbb P (3,4,5,8,11)  &  2\times A_3, A_{10} & 6 \\  \hline
\mathrm{No.}73 & X_{14,18}\subset \mathbb P (2,3,7,9,11)  &  2\times A_2, A_{10} & 8 \\  \hline
\mathrm{No.}74 & X_{14,18}\subset \mathbb P (2,6,7,8,9)  &  5\times A_1, A_2, A_7 & 8\\  \hline
\mathrm{No.}75 & X_{12,20}\subset \mathbb P (4,5,6,7,10)  &  2\times A_1, 2\times A_4, A_6 & 6\\  \hline
\mathrm{No.}76 & X_{16,18}\subset \mathbb P (1,6,8,9,10)  &  A_1, A_2, A_9  & 10 \\  \hline
\mathrm{No.}77 & X_{16,18}\subset \mathbb P (4,6,7,8,9)  &  2\times A_1, 2\times A_3, A_2, A_6 & 6\\  \hline
\mathrm{No.}78 & X_{18,20}\subset \mathbb P (4,5,6,9,14)  &  2\times A_1, A_2, A_{13} & 5\\  \hline
\mathrm{No.}79 & X_{18,20}\subset \mathbb P (4,5,7,9,13)  &  A_6, A_{12} & 4 \\  \hline
\mathrm{No.}80 & X_{18,20}\subset \mathbb P (5,6,7,9,11)  &  A_2, A_6, A_{10} & 4 \\  \hline
\mathrm{No.}81 & X_{18,22}\subset \mathbb P (2,5,9,11,13)  &  A_4, A_{12} & 6 \\  \hline
\mathrm{No.}82 & X_{20,21}\subset \mathbb P (3,4,7,10,17)  &  A_1, A_{16} & 5 \\  \hline
\mathrm{No.}83 & X_{18,30}\subset \mathbb P (6,8,9,10,15)  &  2\times A_1, 2 \times A_2, A_7, A_4 & 5 \\  \hline
\mathrm{No.}84 & X_{24,30}\subset \mathbb P (8,9,10,12,15)  &  A_1, A_3, A_8, A_2, A_4  & 6 \\  \hline
\end{array}
$$
}




\end{document}